# ONSAGER'S CONJECTURE ON THE ENERGY CONSERVATION FOR SOLUTIONS OF EULER EQUATIONS IN BOUNDED DOMAINS

QUOC-HUNG NGUYEN AND PHUOC-TAI NGUYEN

ABSTRACT. The Onsager's conjecture has two parts: conservation of energy, if the exponent is larger than 1/3 and the possibility of dissipative Euler solutions, if the exponent is less or equal than 1/3. The paper proves half of the conjecture, the conservation part, in bounded domains.

*Key words:* Onsager's conjecture; energy conservation; Euler equation.

*MSC:* 35Q31, 76B03.

## 1. INTRODUCTION

The Onsager's conjecture has two parts: conservation of energy, if the exponent is larger than 1/3 and the possibility of dissipative Euler solutions, if the exponent is less or equal than 1/3 (see [13]). In this note, we prove half of the conjecture, the conservation part, for the incompressible Euler equation of the form

$$\begin{cases} \partial_t u + \nabla \cdot (u \otimes u) + \nabla p = 0 & \text{in } \Omega \times (0, T), \\ \nabla \cdot u = 0 & \text{in } \Omega \times (0, T), \\ u(x,t) \cdot n(x) = 0 & \text{on } \partial\Omega \times (0, T), \end{cases} \quad (1.1)$$

where $T > 0$, $\Omega$ is a bounded, connected domain in $\mathbb{R}^d$ ($d \geq 2$) with $C^2$-boundary and $n(x)$ is the outward unit normal vector field to the boundary $\partial\Omega$.

In the absence of a physical boundary (namely the case of whole space $\mathbb{R}^d$ or the case of periodic boundary conditions in the torus $\mathbb{T}^m$), Eyink [11] and Constantin, E and Titi [5] proved that any solution of the incompressible Euler equations must conserve the global kinetic energy if it is Hölder continuous with exponent greater than 1/3.

The optimality of the exponent 1/3 (in the sense that there is a Hölder continuous solution with exponent smaller than 1/3 which dissipates the total kinetic energy) was first established in $\mathbb{R}^2 \times \mathbb{R}$ in the groundbreaking paper of Scheffer [14]. Afterwards, in [15], Shnirelman constructed a compactly supported nontrivial weak solution in $\mathbb{T}^2 \times \mathbb{R}$. Much development in this direction was achieved in a series of celebrated papers of De Lellis and Székelyhidi [6, 7, 8] (see also [4, 12, 3] and references therein).

The first work dealing with the conjecture on the energy conservation for domains with boundaries was recently contributed by Bardos and Titi [1], where the energy conservation was obtained under an assumption on Hölder continuity up to the boundary with exponent greater than 1/3. This assumption is relaxed in [2] where only interior Hölder regularity and continuity of the normal component of the energy flux near the boundary are required. Independently, Drivas and Nguyen [10], by using the Besov spaces, demonstrated the energy conservation which extends the results in [5] and refines the results in [1].

Motivated by the above work, the aim of the present note is to establish the energy conservation under weaker assumptions. Our approach is based on averaging over the strip near the boundary and hence enables us to obtain relaxed conditions.

Before stating the main result, let us introduce the definition of weak solutions of (1.1).

We say that $(u, p)$ is a *weak solution* of (1.1) in $\Omega \times (0, T)$ if $u \in C_w((0, T), H(\Omega))$, $p \in L^1_{loc}(\Omega \times (0, T))$ and

$$\int_0^T \int_\Omega (u \cdot \partial_t \varphi + u \otimes u : \nabla \varphi + p \nabla \cdot \varphi) dx dt = 0 \qquad (1.2)$$

for all test vector fields $\varphi \in C_0^\infty(\Omega \times (0, T))$. Here $H(\Omega)$ is the completion in $L^2(\Omega)$ of the space

$$\{w \in C_c^\infty(\Omega, \mathbb{R}^d) : \nabla \cdot w = 0\}.$$

In the sequel, $d(x)$ is the distance from $x$ to $\partial \Omega$, $\Omega_r := \{x \in \Omega : d(x) > r\}$ for any $r \geq 0$ and $\fint_E f dx := \frac{1}{\mathcal{L}^d(E)} \int_E f dx$ for any Borel set $E \subset \mathbb{R}^d$.

Since $\Omega$ is a bounded, connected domain with $C^2$-boundary, so we find $r_0 > 0$ and a unique $C_b^1$-vector function $n : \Omega \backslash \Omega_{r_0} \to S^{d-1}$ such that the following holds true: for any $r \in [0, r_0)$, $x \in \Omega_r \backslash \Omega_{r_0}$ there exists a unique $x_r \in \partial \Omega_r$ such that $d(x, \partial \Omega_r) = |x - x_r|$ and $n(x)$ is the outward unit normal vector field to the boundary $\partial \Omega_r$ at $x_r$. In particular, for any $0 < r_1 < r_2 < r_0$, and $\varphi \in L^1(\Omega_{r_1} \backslash \Omega_{r_2})$, by the Coarea formula we have

$$\int_{\Omega_{r_1} \backslash \Omega_{r_2}} \varphi(x) dx = \int_{r_1}^{r_2} \int_{\partial \Omega_\tau} \varphi(\theta) d\mathcal{H}^{d-1}(\theta) d\tau. \qquad (1.3)$$

Our main theorem is the following:

**Theorem 1.1.** *Let $(u, p)$ is a weak solution of problem (1.1). Assume that for any $\varepsilon > 0$, there exist $\sigma_\varepsilon > 1/3$ and $C_\varepsilon > 0$ such that*

$$\sup_{|h| < \varepsilon/2} |h|^{-\sigma_\varepsilon} \|u(\cdot + h, \cdot) - u(\cdot, \cdot)\|_{L^3(\Omega_\varepsilon \times (0,T))} < C_\varepsilon < \infty. \qquad (1.4)$$

*In addition, assume that for some $\varepsilon_0 \in (0, r_0)$ small, $u \in L^3((\Omega \backslash \Omega_{\varepsilon_0}) \times (0, T))$, $p \in L^{3/2}((\Omega \backslash \Omega_{\varepsilon_0}) \times (0, T))$ and*

$$\left(\int_0^T \fint_{\Omega \backslash \Omega_\varepsilon} |u(x,t)|^3 dx dt\right)^{\frac{2}{3}} \left(\int_0^T \fint_{\Omega \backslash \Omega_\varepsilon} |u(x,t) \cdot n(x)|^3 dx dt\right)^{\frac{1}{3}} = \circ(1), \qquad (1.5)$$

$$\left(\int_0^T \fint_{\Omega \backslash \Omega_\varepsilon} |p(x,t)|^{\frac{3}{2}} dx dt\right)^{\frac{2}{3}} \left(\int_0^T \fint_{\Omega \backslash \Omega_\varepsilon} |u(x,t) \cdot n(x)|^3 dx dt\right)^{\frac{1}{3}} = \circ(1), \qquad (1.6)$$

*as $\varepsilon \to 0$.*

*Then the energy is globally conserved, namely*

$$\int_\Omega |u(x,t)|^2 dx = \int_\Omega |u(x,0)|^2 dx \quad \forall t \in [0, T).$$

**Remark 1.2.** *(i) Conditions (1.5) and (1.6) can be replaced by the following conditions*

$$\limsup_{\varepsilon \to 0} \int_0^T \fint_{\Omega \backslash \Omega_\varepsilon} |u(x,t)|^3 dx dt + \limsup_{\varepsilon \to 0} \int_0^T \fint_{\Omega \backslash \Omega_\varepsilon} |p(x,t)|^{\frac{3}{2}} dx dt < \infty$$

and
$$\liminf_{\varepsilon \to 0} \int_0^T \fint_{\Omega \setminus \Omega_\varepsilon} |u(x,t) \cdot x|^3 dx dt = 0.$$

(ii) Put $\tilde{u}(x,t) = u(x,t) \cdot n(x)$ for any $x \in \Omega \setminus \Omega_{r_0}$. If the function
$$h : \varepsilon \mapsto h(\varepsilon) = \|\tilde{u}\|_{L^3((\Omega \setminus \Omega_\varepsilon) \times (0,T))}$$
satisfies $h(\varepsilon) \leq C\varepsilon^{2/3}$ for every $\varepsilon \in (0, r_0)$ with some $C > 0$, then conditions (1.5) and (1.6) are fulfilled.

(iii) Our conditions are weaker than those in [10]. In particular, in Theorem 1.1, the velocity field $u$ and the pressure $p$ are required to be in $L^3$ and $L^{3/2}$ in a strip near the boundary, while in [10, Theorem 1] they are required to be bounded in a strip near the boundary. Moreover, if $u \in L^3(0,T, L^\infty(\Omega \setminus \Omega_{\varepsilon_0}))$ and $p \in L^{3/2}(0,T, L^\infty(\Omega \setminus \Omega_{\varepsilon_0}))$ then our conditions (1.5) and (1.6) can be deduced from the continuity of the wall-normal velocity at the boundary in [10, Theorem 1].

## 2. Proof of Theorem 1.1

For the sake of simplicity, we proceed as if the solution is differentiable in time. The extra argument needed to mollify in time can be found in [9]. Let $\{\varrho_\varepsilon\}$ be a standard sequence of mollifiers in $\mathbb{R}^d$ satisfying $\mathrm{supp}\varrho_\varepsilon \subset B_\varepsilon$ where $B_\varepsilon$ is the ball of radius $\varepsilon$ and center at 0. For any function $f \in L^1_{loc}(\Omega \times (0,T))$, we set for any $x \in \Omega_\varepsilon$,

$$f^\varepsilon(x,t) := (\varrho_\varepsilon \star f(\cdot,t))(x) = \int_\Omega \varrho_\varepsilon(x-y)f(y,t)dy = \int_{\mathbb{R}^d} \varrho_\varepsilon(x-y)f(y,t)dy. \quad (2.1)$$

From (1.1), we derive
$$\partial_t u^\varepsilon + \nabla \cdot (u \otimes u)^\varepsilon + \nabla p^\varepsilon = 0, \ \nabla \cdot u^\varepsilon = 0 \quad \text{in } \Omega_{2\varepsilon} \times (0,T).$$

This implies,
$$\partial_t u^\varepsilon + \nabla \cdot (u^\varepsilon \otimes u^\varepsilon) + \nabla p^\varepsilon + \nabla \cdot r_\varepsilon = 0 \quad \text{in } \Omega_{2\varepsilon} \times (0,T). \quad (2.2)$$

Here we have used the fact (see [5]) that
$$(u \otimes u)^\varepsilon = u^\varepsilon \otimes u^\varepsilon + r_\varepsilon \quad \text{in } \Omega_{2\varepsilon} \times (0,T),$$
with
$$r_\varepsilon(x,t) := \int_\Omega \varrho_\varepsilon(y)(u(x-y,t) - u(x,t)) \otimes (u(x-y,t) - u(x,t))dy$$
$$- (u - u^\varepsilon)(x,t) \otimes (u - u^\varepsilon)(x,t), \quad (x,t) \in \Omega_{2\varepsilon} \times (0,T).$$

Thus, it follows from (2.2) that
$$\partial_t \left(\frac{1}{2}|u^\varepsilon|^2\right) + \nabla \cdot F_\varepsilon = G_\varepsilon \quad \text{in } \Omega_{2\varepsilon} \times (0,T). \quad (2.3)$$

where
$$F_\varepsilon := \left(\frac{1}{2}|u^\varepsilon|^2 + p^\varepsilon\right)u^\varepsilon + u^\varepsilon \cdot r_\varepsilon \quad \text{and} \quad G_\varepsilon := \nabla u^\varepsilon : r_\varepsilon.$$

Let $0 < \varepsilon < \varepsilon_1/10 < \varepsilon_2/10 < r_0/100$. Integrating (2.3) over $\Omega_{\varepsilon_2} \times (0,t)$ yields

$$\frac{1}{2}\int_{\Omega_{\varepsilon_2}} |u^\varepsilon(x,t)|^2 dx - \frac{1}{2}\int_{\Omega_{\varepsilon_2}} |u^\varepsilon(x,0)|^2 dx$$
$$= -\int_0^t \int_{\partial\Omega_{\varepsilon_2}} F_\varepsilon(\theta,s) \cdot n(\theta) d\mathcal{H}^{d-1}(\theta) ds + \int_0^t \int_{\Omega_{\varepsilon_2}} G_\varepsilon(x,s) dx ds.$$

So, for any $\varepsilon_3 > 0$ small, by integrating over $(\varepsilon_1, \varepsilon_1 + \varepsilon_3)$ in $\varepsilon_2$ and using (1.3) for the first term on the right-hand side, we obtain

$$\frac{1}{2\varepsilon_3}\int_{\varepsilon_1}^{\varepsilon_1+\varepsilon_3}\int_{\Omega_{\varepsilon_2}} |u^\varepsilon(x,t)|^2 dx d\varepsilon_2 - \frac{1}{2\varepsilon_3}\int_{\varepsilon_1}^{\varepsilon_1+\varepsilon_3}\int_{\Omega_{\varepsilon_2}} |u^\varepsilon(x,0)|^2 dx d\varepsilon_2$$
$$= -\int_0^t \frac{1}{\varepsilon_3}\int_{\Omega_{\varepsilon_1}\setminus\Omega_{\varepsilon_1+\varepsilon_3}} F_\varepsilon(x,s) \cdot n(x) dx ds + \int_0^t \frac{1}{\varepsilon_3}\int_{\varepsilon_1}^{\varepsilon_1+\varepsilon_3}\int_{\Omega_{\varepsilon_2}} G_\varepsilon(x,s) dx d\varepsilon_2 ds. \tag{2.4}$$

Thanks to (1.4), we derive

$$\left|\int_0^t \int_{\Omega_{\varepsilon_1}\setminus\Omega_{\varepsilon_1+\varepsilon_3}} u^\varepsilon(x,s) \cdot r_\varepsilon(x,s) \cdot n(x) dx ds\right|$$
$$\leq C_{\varepsilon_1,\varepsilon_3}\varepsilon^{2\sigma_{\varepsilon_1,\varepsilon_3}}\|u^\varepsilon\|_{L^3((0,T)\times(\Omega_{\varepsilon_1}\setminus\Omega_{\varepsilon_1+\varepsilon_3}))} \tag{2.5}$$
$$\leq C_{\varepsilon_1,\varepsilon_3}\varepsilon^{2\sigma_{\varepsilon_1,\varepsilon_3}}\|u\|_{L^3((0,T)\times(\Omega_{\varepsilon_1}\setminus\Omega_{\varepsilon_1+\varepsilon_3}))},$$

for some constant $C_{\varepsilon_1,\varepsilon_3} > 0$ and $\sigma_{\varepsilon_1,\varepsilon_3} > \frac{1}{3}$.

Similarly, since $\varepsilon_1 < \varepsilon_2$ and by (1.4),

$$\int_0^t \int_{\varepsilon_1}^{\varepsilon_1+\varepsilon_3}\int_{\Omega_{\varepsilon_2}} |G_\varepsilon(x,s)| dx d\varepsilon_2 ds \leq \varepsilon_3 \int_0^t \int_{\Omega_{\varepsilon_1}} |G_\varepsilon(x,s)| dx ds \leq \varepsilon_3 C_{\varepsilon_1}\varepsilon^{3\sigma_{\varepsilon_1}-1}, \tag{2.6}$$

for some constant $C_{\varepsilon_1} > 0$ and $\sigma_{\varepsilon_1} > \frac{1}{3}$.

Letting sucessively $\varepsilon \to 0$ and then $\varepsilon_1 \to 0$ in (2.4) and taking account into (2.5) and (2.6), we obtain

$$\frac{1}{2\varepsilon_3}\int_0^{\varepsilon_3}\int_{\Omega_{\varepsilon_2}} |u(x,t)|^2 dx d\varepsilon_2 - \frac{1}{2\varepsilon_3}\int_0^{\varepsilon_3}\int_{\Omega_{\varepsilon_2}} |u(x,0)|^2 dx d\varepsilon_2$$
$$= -\int_0^t \frac{1}{\varepsilon_3}\int_{\Omega\setminus\Omega_{\varepsilon_3}} \frac{1}{2}|u(x,s)|^2 u(x,s) \cdot n(x) dx ds - \int_0^t \frac{1}{\varepsilon_3}\int_{\Omega\setminus\Omega_{\varepsilon_3}} p(x,s) u(x,s) \cdot n(x) dx ds.$$

Using Holder's inequality for the terms on the right-hand side yields

$$\left|\frac{1}{2\varepsilon_3}\int_0^{\varepsilon_3}\int_{\Omega_{\varepsilon_2}} |u(x,t)|^2 dx d\varepsilon_2 - \frac{1}{2\varepsilon_3}\int_0^{\varepsilon_3}\int_{\Omega_{\varepsilon_2}} |u(x,0)|^2 dx d\varepsilon_2\right|$$
$$\leq C\left(\int_0^T \fint_{\Omega\setminus\Omega_{\varepsilon_3}} |u(x,t)|^3 dx dt\right)^{\frac{2}{3}} \left(\int_0^T \fint_{\Omega\setminus\Omega_{\varepsilon_3}} |u(x,t) \cdot n(x)|^3 dx dt\right)^{\frac{1}{3}} \tag{2.7}$$
$$+ C\left(\int_0^T \fint_{\Omega\setminus\Omega_{\varepsilon_3}} |p(x,t)|^{\frac{3}{2}} dx dt\right)^{\frac{2}{3}} \left(\int_0^T \fint_{\Omega\setminus\Omega_{\varepsilon_3}} |u(x,t) \cdot n(x)|^3 dx dt\right)^{\frac{1}{3}}.$$

Here we have used the fact that $\mathcal{L}^d(\Omega\setminus\Omega_{\varepsilon_3}) \simeq \varepsilon_3$. Owing to (1.5) and (1.6), by letting $\varepsilon_3 \to 0$ in (2.7), we conclude

$$\int_\Omega |u(x,t)|^2 dx = \int_\Omega |u(x,0)|^2 dx \quad \forall t \in [0,T).$$

The proof is complete.

**Acknowledgements.** The authors are grateful to Emil Wiedemann for helpful comments. We also would like to thank the anonymous referee for constructive comments which help to improve the note remarkably.

## References

[1] Bardos, C., and Titi, E.: Onsager's Conjecture for the Incompressible Euler Equations in Bounded Domains. Archive for Rational Mechanics and Analysis 228.1, 197-207 (2018).
[2] Bardos, C., Titi, E., Wiedemann, E.: Onsager's Conjecture with Physical Boundaries and an Application to the Viscosity Limit, preprint arXiv:1803.04939, (2018).
[3] Buckmaster,T., De Lellis, C., Székelyhidi Jr. L. and Vicol,V.: Onsagers conjecture for admissible weak solutions, arXiv:1701.08678 (2017).
[4] Buckmaster, T., De Lellis, C., Isett, P. and Székelyhidi, Jr. L.: Anomalous dissipation for 1/5-Hölder Euler flows. Ann. of Math. (2) 182, no. 1, 127-172 (2015).
[5] Constantin, P., W. E, and Titi, E.: Onsager's conjecture on the energy conservation for solutions of Euler's equation. Comm. Math. Phys. 165: 207-209 (1994).
[6] De Lellis, C. and Székelyhidi Jr., L.: The h-principle and the equations of fluid dynamics, B. Am. Math. Soc. 49: 347-375 (2012).
[7] De Lellis, C. and Székelyhidi, Jr. L.: Dissipative continuous Euler flows. Invent. Math. 193, no. 2, 377-407 (2013).
[8] De Lellis, C. and Székelyhidi, Jr. L.: Dissipative Euler flows and Onsager's conjecture. J. Eur. Math. Soc. (JEMS) 16, no. 7, 1467-1505 (2014).
[9] Drivas, T.D. and Eyink, G.L.: An Onsager singularity theorem for turbulent solutions of compressible Euler equations. Comm. Math. Phys. 1-31 (2017).
[10] Drivas, T.D. and Nguyen, H. Q.: Onsager's conjecture and anomalous dissipation on domains with boundary, preprint arXiv: 1803.05416v2 (2018).
[11] Eyink, G.L.: Energy dissipation without viscosity in ideal hydrodynamics I. Fourier analysis and local energy transfer, Phys. D 78, 222-240 (1994).
[12] Isett, P.: A proof of Onsager's conjecture, preprint arXiv:1608.08301, (2016).
[13] Onsager, L.: Statistical Hydrodynamics. Nuovo Cimento (Supplemento), 6, 279 (1949).
[14] Scheffer, V.: An inviscid flow with compact support in space-time. J. Geom. Anal. 3(4), 343-401 (1993).
[15] Shnirelman, A.: On the nonuniqueness of weak solution of the Euler equation. Comm. Pure Appl. Math. 50(12), 1261-1286 (1997).

Scuola Normale Superiore, Piazza dei Cavalieri 7, I-56100 Pisa, Italy
*E-mail address*: quochung.nguyen@sns.it

Department of Mathematics, Masaryk University, Brno, Czech Republic
*E-mail address*: ptnguyen@math.muni.cz;   nguyenphuoctai.hcmup@gmail.com